\documentclass[a4paper]{amsart}

\usepackage {amssymb}
\usepackage{verbatim}
\usepackage{amscd}

\theoremstyle{plain}

\theoremstyle{definition}

\def\deg{\operatorname{deg}}
\def\det{\operatorname{det}}

\def\Hom{\operatorname{Hom}}

\def\Ind{\operatorname{Ind}}

\def\deg{\operatorname{deg}}

\def\symm{\operatorname{symm}}
\def\Res{\operatorname{Res}}

\newcommand{\field}{\mathbb}

\newcommand{\Z}{{\field Z}}

\newcommand{\Spec}{{\mathrm{Spec}}}

\newcommand{\fra}{\mathfrak{a}}
\newcommand{\frb}{\mathfrak{b}}
\newcommand{\frc}{\mathfrak{c}}

\newcommand{\fre}{\mathfrak{e}}
\newcommand{\frf}{\mathfrak{f}}
\newcommand{\frg}{\mathfrak{g}}
\newcommand{\frh}{\mathfrak{h}}

\newcommand{\frk}{\mathfrak{k}}
\newcommand{\frl}{\mathfrak{l}}
\newcommand{\frm}{\mathfrak{m}}

\newcommand{\fro}{\mathfrak{o}}
\newcommand{\frp}{\mathfrak{p}}
\newcommand{\frq}{\mathfrak{q}}

\newcommand{\frs}{\mathfrak{s}}

\newcommand{\fru}{\mathfrak{u}}

\newcommand{\bbC}{\mathbb{C}}
\newcommand{\bbD}{\mathbb{D}}

\newcommand{\bbR}{\mathbb{R}}

\newcommand{\bbZ}{\mathbb{Z}}

\newcommand{\caD}{\mathcal{D}}
\newcommand{\caE}{\mathcal{E}}

\newcommand{\caH}{\mathcal{H}}

\newcommand{\caO}{\mathcal{O}}

\newcommand{\caS}{\mathcal{S}}

\newcommand{\caW}{\mathcal{W}}

\newcommand{\del}{\partial}

\begin{document}

\title[Erratum and addendum to: Invariant Differential Operators]{Erratum and Addendum to: Invariant Differential Operators and Eigenspace Representations on an Affine Symmetric Space}

\author{Jing-Song Huang}
\address{Department of Mathematics, Hong Kong University of Science and Technology,
Clear Water Bay, Kowloon, Hong Kong SAR, China}
\email{mahuang\char'100ust.hk}



\thanks{I thank heartily Nolan Wallach for helpful discussions and the referees for their useful suggestions.
}

\maketitle


The purpose of this erratum and addendum is to correct the errors in \cite{H}.  It consists of five components:

1. Lemma 7.1 and Proposition 7.2 are wrong and discarded;

2. A new proof of existence $\lambda(\xi)$ in (7.1) without Proposition 7.2;

3. Definition of a new bijection in Theorem 5.2 and a proof by a new technique;

4. A new proof of Theorem 5.5 based on the new bijection in Theorem 5.2;

5. Correction to the list of exceptional simple pairs in Proposition 3.1.

\noindent The main results of [1] remain true as stated.  We also add a final remark on generalization.

\subsection*{1. Discarding Lemma 7.1 and Proposition 7.2.}
Lemma 7.1 and Proposition 7.2 missed a necessary condition, namely the assumption that $\nu_0$ is invariant
under the induced action of $N_K(\fra_{\frq_\bbC})$.  In case $\nu_0=\rho(\frm)$, the lemma
is Proposition 3.1 of [2].
However, this condition fails for general $\nu_0$.  We discard both Lemma 7.1 and Proposition 7.2.

\subsection*{2. A Proof of Existence of $\lambda(\xi)$ in (7.1)}
Recall that $\fra_\frq$ is a maximal abelian subalgebra in $\frp\cap \frq$ and $\frc\supseteq \fra_\frq$ is a Cartan subalgebra of $\frg$.
Let $\frc_0$ denote the orthogonal complement of $\fra_\frq$ in $\frc$.
Then $\frc=\frc_0\oplus\fra_\frq$ and $\frc_\bbC={\frc_0}_\bbC\oplus {\fra_\frq}_\bbC$. Let $t=(t_1,\ldots,t_{n-r})\in{\frc_0}^*_\bbC$ and $x=(x_1,\ldots,x_r)\in{\fra_\frq}^*_\bbC$ denote the coordinates.

We regard symmetric algebra $S(\frc)$ (resp. $S(\fra_\frq)$) as algebra of complex valued polynomial functions on $\frc_\bbC^*$
(resp.  ${\fra_\frq}^*_\bbC$).
Let $I(\frc)=S(\frc)^{W(\frc)}$ and $I(\fra_\frq)=S(\fra_\frq)^{W(\fra_\frq)}$ be the Weyl group invariants. Suppose that $U_1,\cdots,U_r\in I(\frc)$ are homogeneous polynomials in $(t;x)\in \frc^*_\bbC$
such that the restriction to ${\fra_\frq}_\bbC^*$
$$W_1(x)=\Res^\frc_{\fra_\frq}U_1=U_1(0;x),\cdots,W_r(x)=\Res^\frc_{\fra_\frq}U_r=U_r(0;x)$$ are algebraically independent and
$I(\fra_\frq)$ is a finite free module over $\bbC[W_1,\cdots,W_r]$ of rank $d$.
Consider the Jacobian determinant
$$J(t;x)=\det\left[ \frac{\del U_i(t;x)}{\del x_j} \right].$$
Regarding $U_i(t;x)=\sum a_{j_1,\dots,j_{n-r}}(x)t_1^{j_1}\cdots t_{n-r}^{j_{n-r}}$ as a polynomial in $t$, the constant
term is $W_i(x)=U_i(0;x)$. The
coefficients for non-constant terms in $t$ are polynomials in $x$ with strictly smaller degrees.
It follows that $\frac{\del U_i(t;x)}{\del x_j}$ is a homogeneous
polynomial in $(t;x)$ of degree $\deg U_i-1$.
Thus, $J(t;x)$ is a homogeneous polynomial in $(t;x)$ of degree equal to $(\deg U_1-1) \cdots(\deg U_r -1)$.
By algebraically independence of $W_1,\cdots, W_r$, we have that
$$\det\left[ \frac{\del W_i(x)}{\del x_j} \right]=\det\left[ \frac{\del U_i(t;x)}{\del x_j} \right]_{t=0}=J(0;x)$$
is a nonzero homogeneous polynomial in $x$.

The expansion of $J(t;x)$
in terms of polynomial in $t$ has the constant term $J(0;x)$
and the  coefficients for other terms are polynomials in $x$ with strictly less degrees.
Thus, for a fixed $\zeta\in{\frc_0}^*_\bbC$,  $J(\zeta;x)$ is a nonzero polynomial with the leading term $J(0;x)$.
In particular,  $U_1(\zeta;x),\ldots,U_r(\zeta;x)$ are algebraically independent. Set $A_\zeta=\bbC[U_1(\zeta;x),\ldots,U_r(\zeta;x)]$.
Then the map
$$\Gamma\colon A_\zeta\rightarrow S(\fra_\frq),  \ \
f\mapsto f(U_1(\zeta;x),\ldots,U_r(\zeta;x))$$ is an injective algebra homomorphism.
We claim that $S(\fra_\frq)$ is integral over $A_\zeta$. It is enough to show any
polynomial $P(x)\in S(\fra_\frq)$ is integral over $A_\zeta$.  We prove this by induction on
$\deg P(x)$.  If $\deg P(x)=0$, this is obvious.  Assume that any polynomial $P(x)$ with  $\deg P\leq m$ is integral over $A_\zeta$.  Now we look at a $P(x)$ with $\deg P(x)=m+1$.
Since $S(\fra_\frq)$ is integral over $I(\fra_\frq)$, we may reduce to the case that
$P(x)$ in $I(\fra_\frq)$.  Since $I(\fra_\frq)$ is integral over $A_0=\bbC[W_1,\ldots, W_r]$,
we may further reduce to the case $P(x)$ is in $A_0$.  More precisely, suppose that
$e_1,\ldots,e_m$ are generators of $S(\fra_\frq)$ as a $A_0$-module. Then there are
$P_i(x)=F_i(W_1,\ldots,W_r)\in A_0$ for some homogeneous
polynomial $F$ such that degree of $P_i(x)\leq m+1$ and
$$P(x)=P_1(x)e_1+\cdots +P_m(x)e_m.$$
It follows that for each $i\in\{1,\ldots,m\}$ $$Q_i(x)=F_i(U_1(\zeta;x),\ldots,U_r(\zeta;x))-P_i(x)=F_i(U_1(\zeta;x),\ldots,U_r(\zeta;x))-F_i(U_1(0;x),\ldots,U_r(0;x))$$
has degree at most $m$.  By induction hypothesis, $Q_i(x)$ is integral over $A_\zeta$,
hence $P_i(x)$ is integral over $A_\zeta$ for $i=1,\ldots, m$.  Then $P(x)$ is integral over $A_\zeta$-module generated by
$P_1,\ldots, P_m,e_1,\ldots , e_m$
and hence it is integral over $A_\zeta$.

The injective algebra homomorphism $\Gamma: A_\zeta \rightarrow S(\fra_\frq)$ induces a finite map
$\phi: \Spec S(\fra_\frq) \rightarrow \Spec A_\zeta$.   By identifying ${\fra_\frq}^*_\bbC$
with $\bbC^r$ we have
$$\phi\colon \bbC^r\rightarrow \bbC^r,  \ \ x\mapsto (U_1(\zeta;x),\ldots,U_r(\zeta;x)).$$
The integrality of $S(\fra_\frq)$ over $A_\zeta$ implies that $\phi$ is surjective.
More precisely, let $a=(a_1,\ldots,a_r)\in \bbC^r$ determine the maximal ideal
$$I_a=(U_1(\zeta;x)-a_1,\ldots ,U_r(\zeta;x)-a_r)\bbC[U_1(\zeta;x),\ldots,U_r(\zeta;x)]$$
in $A_\zeta$.  Then there exists a maximal ideal $I_b$ in $\bbC[x_1,\ldots,x_r]$  with $b\in \bbC^r$ and
$$I_b=(x_1-b_1,\ldots ,x_r-b_r)\bbC[x_1,\ldots,x_r]$$
such that $I_b\cap \Gamma(\bbC[U_1(\zeta;x),\ldots, U_r(\zeta;x)])=\Gamma(I_a)$.
In other words, we have
$$U_1(\zeta,b_1)=a_1,\ldots,U_r(\zeta,b_r)=a_r.$$

Let $\gamma\colon Z(\frg)\rightarrow I(\frc)$ be the Harish-Chandra isomorphism.
Let ${\bf D}=\bbC[D_1,\ldots ,D_r]$ be the polynomial subalgebra of $Z(\frg)$
generated by $D_i$ with $\gamma(D_i)=U_i$ ($i=1,\cdots, r$).
For  $\lambda \in {\fra_\frq}_\bbC^*$,
we define the character
$$\chi_\lambda\colon {\bf D}=\bbC[D_1,\cdots , D_r]\rightarrow \bbC \text{ by }\chi_\lambda(D_i)=\langle\gamma(D_i),\lambda\rangle=U_i(0;\lambda).$$ For a given
finite-dimensional unitary representation $\xi\in  \widehat{M}_{fu}$, denote by
$\Lambda_\xi\in {\frc_0}^*_\bbC$ its infinitesimal character.  Then there exists $\lambda(\xi)\in {\fra_\frq}^*_\bbC$ such that
$$U_i(\Lambda_\xi;\lambda(\xi))=U_i(0;\lambda), i=1,\cdots, r.$$  Thus,  $\gamma(D_i,\Lambda_\xi+\lambda(\xi))=\chi_\lambda(D_i)$
 for each $i$, and it follows that
$$\gamma(Z,\Lambda_\xi+\lambda(\xi))=\chi_\lambda(Z), \text{\ for all } Z\in \bbC[D_1,\ldots, D_r].$$
This establishes (7.1).

For any $\xi\in  \widehat{M}_{fu}$, the polynomial equation
$$ J(\Lambda_\xi;x)=\det\left[ \frac{\del U_i(\Lambda_\xi;x)}{\del x_j} \right]=0$$ 
defines a hypersurface in
${\fra_\frq}^*_\bbC$.  We say that $x\in {\fra_\frq}^*_\bbC$ is {\it unramified} if
$J(\Lambda_\xi;x)\neq 0$ for all $\xi\in \widehat{M}_{fu}$.  The ramified points (the complement to
unramified points) are collection of locally finite hypersurfaces defined by $J(\Lambda_\xi;x)=0$, since
any ideal of $S(\fra_\frq)$ is finitely generated and therefore each point $x\in {\fra_\frq}^*_\bbC$ can
lay in only finite many of such hypersurfaces.
We say $x\in {\fra_\frq}^*_\bbC$ is {\it generic} if it is unramified and in addition it is not on the hyperplanes
defined by $\langle x,\alpha \rangle \in \bbZ$.
If $x$ is unramified, then the set of fibre $\phi^{-1}(y)$ with $y=\phi(x)$ has exactly $|W(\fra_\frq)|\cdot d$ elements.
This follows from the fact that $J(\xi;x)\neq 0$ for all $x\in \phi^{-1}(y)$
and hence $\phi$ maps a neighborhood of $x$ homeomorphically
to a neighborhood of $y=\phi(x)$.  Since each $W(\fra_\frq)$-orbits can have at most $|W(\fra_\frq)|$ points, there are at least
$d$ points $\lambda(\xi)_1,\ldots \lambda(\xi)_d\in {\fra_\frq}^*_\bbC$
in distinct $W(\fra_\frq)$-orbits, such that
$$\gamma(Z,\Lambda_\xi+\lambda(\xi)_i)=\chi_\lambda(Z), \text{\ for all } Z\in \bbC[D_1,\ldots, D_r].$$

The character $\chi_\lambda$ of $\bf D$ ($\lambda\in {\fra_\frq}^*_\bbC$) or simply
$\lambda$ is said to be generic if all points  $\phi^{-1}(\lambda)$ are generic. 
If $\lambda$ is generic, then the principal series $\pi_{\xi,\lambda(\xi)_i}$ ($i=1,\ldots ,d$)
are irreducible and non-isomorphic.  We also note that
if  $\nu_0$ is unramified, then the function $\nu\mapsto \lambda=(U_1(\Lambda_\xi;\nu),\ldots,U_r(\Lambda_\xi;\nu))$
for  $\nu\in {\fra_\frq}^*_\bbC$
has an inverse $\psi$ in some neighborhood of $\nu_0$
by the inverse function theorem, and
the inverse $\psi\colon\lambda \mapsto \nu=\lambda(\xi)$ is holomorphic in some neighborhood of
$\lambda_0=(U_1(\Lambda_\xi,\nu_0), \ldots, U_r(\Lambda_\xi,\nu_0))$.

\subsection*{3. The bijection in Theorem 5.2.}

I am grateful to E. van den Ban and P. Delorme for pointing out an error in determining the
$\tau$-radical part of a differential operator in Section 5 \cite{H}.  The error is
the claim that `$\bar D$ is right invariant under $G_+$' in 5th line on Page 719.
The whole paragraph
(from Line -12 Page 718 to Line 13 page 719) containing this wrong claim is discarded.
We define a new bijection for Theorem 5.2 by restriction of Taylor expansions.

In order to define this new bijection we set up two linear maps $\Gamma_\frg$ and $\Gamma_{\frg_+}$ regarding $U(\frg)$ and $U(\frg_+)$ respectively.
Recall from (3.6) that the restriction map
$p\colon S(\frc)^{W(\frc)}\rightarrow S(\fra_{\frq})^{W(\fra_\frq)}$ is the
composition of the following two maps
$$ I(\frc)=S(\frc)^{W(\frc)}\rightarrow S(\frb)^{W(\frb)}\rightarrow S(\fra_{\frq})^{W(\fra_\frq)}=I(\fra_\frq).$$
This is a surjective map unless $\frg$ contains exceptional simple Lie algebras of type $E_6, E_7$ or $E_8$.
In any case, there are homogeneous polynomials $U_1, \cdots U_r\in I(\frc)$ such that the restriction
$W_i=\Res^\frc_{\fra_\fra}U_i$ ($i=1, \cdots, r$) are algebraically independent and
$I(\fra_\frq)$ is a free module over the subalgebra $\bbC[W_1\ldots, W_r]$ of rank $d$,
where $d$ is 1 except for the exceptional cases which are determined in Proposition 3.1.

We regard symmetric algebra $S(\frg)$ as an algebra of complex polynomial functions on $\frg^*$.
Denote by $\caO(\frg)$ the space of complex value polynomial functions on $\frg$. By using the Killing form,
we identify the symmetric algebra $S(\frg)$ with $\caO(\frg)$,
$S(\frg_+)$ with $\caO(\frg_+)$ and $S(\frp\cap\frq)$ with $\caO(\frp\cap\frq)$, etc.
By Chevalley's restriction theorem, we have the following algebra isomorphisms:
$$\Res^\frg_\frc\colon S(\frg)^G\rightarrow S(\frc)^{W(\frc)} \text{ and }
\Res^{\frp\cap\frq}_{\fra_\frq}\colon S(\frp\cap\frq)^{K\cap H}\rightarrow S(\fra_\frq)^{W_{K\cap H}}.$$
Denote by $\iota(U_i)$ the inverse of $U_i$ under
and $\iota(W_i)$ the inverse of $W_i$,  namely,
$$\Res^\frg_\frc \iota(U_i)=U_i  \text{ and } \Res^{\frp\cap\frq}_{\fra_\frq}\iota(W_i)=W_i.$$
Let $\caD=\bbC[\iota(U_1),\ldots ,\iota(U_r)]$.
Let $E$ be a subspace of $S(\frp\cap\frq)^{K\cap H}$
so that $S(\frp\cap\frq)^{K\cap H}$ is a free $\caD$-module generated by $E$, where the
action of $\iota(U_i)$ is multiplication by $\Res^\frg_{\frp\cap\frq}\iota(U_i)=\iota(W_i)$.
We have $\dim E = |W(\fra_\frq)/W_{K\cap H}|\cdot d = |\caW|\cdot d$.
It follows that the linear map
$$ S(\frk)\otimes \caH(\frp\cap\frq)\otimes E\otimes {\caD}\otimes S(\frh)\rightarrow S(\frg) $$
given by $k\otimes h^+\otimes e \otimes D \otimes h\mapsto k h^+eDh$ is a surjection.
Recall that
${\bf D}=\bbC[D_1,\ldots ,D_r]$ is the polynomial subalgebra of $Z(\frg)$ generated
by $D_i$ with $\gamma(D_i)=U_i$ ($i=1,\cdots, r$).
Let $U^j(\frg)\subset U^{j+1}(\frg)$ be the standard filtration of the universal enveloping algebra.
We compare the grade with the filtration and argue as in [W, 11.2.2], and
conclude that the linear map
$$ \Gamma_\frg\colon U(\frk)\otimes \caH(\frp\cap\frq)\otimes E\otimes {\bf D}\otimes U(\frh)\rightarrow U(\frg) $$
given by $k\otimes h^+\otimes e \otimes D \otimes h\mapsto k \symm(h^+)\symm(e)Dh$ is a surjection.
We remark that in the above linear surjection $\Gamma_\frg$ the algebra $\bf D$ 
can be replaced by the polynomial subalgebra of $Z(\frg)$ generated by
$\symm (\iota(U_i))$, since $\gamma(\symm (\iota(U_i)))$ equals to $U_i$ modulo lower degree terms.

Now we define the second map $\Gamma_+$. Let $\frc_+$ be a Cartan subalgebra of $\frg_+$ containing $\fra_\frq$.
Note that the restriction $p$ is also equal to the
the composition of the following two maps
$$S(\frc)^{W(\frc)}\rightarrow S(\frc_+)^{W(\frc_+)}\rightarrow S(\fra_{\frq})^{W(\fra_\frq)}.$$
Set $U'_i=\Res^{\frc}_{\frc_+} U_i$.  Then $\Res^{\frc_+}_{\fra_\frq}U'_i=W_i$.
By Chevalley's restriction theorem, the restrictions
$$\Res^{\frg_+}_{\frc_+}\colon S(\frg_+)^{G_+}\rightarrow S(\frc_+)^{W(\frc_+)} $$ is
an algebra isomorphism.  Denote by $\iota(U'_i)$ the inverse of $U'_i$, namely,
$\Res^{\frg_+}_{\frc_+} \iota(U'_i)=U'_i$.
It follows that $\Res^{\frg_+}_{\frp\cap\frq}\iota(U'_i)=\iota(W_i)$.
Let $\caD'=\bbC[\iota(U'_1),\ldots ,\iota(U'_r)]$.
Then $S(\frp\cap\frq)^{K\cap H}$ is a free $\caD'$-module generated by $E$ (the
action of $\iota(U'_i)$ is multiplication by $\iota(W_i)$).
Then the map $h^+\otimes e \otimes D \otimes h\mapsto h^+eDh$ defines a linear bijection
$$\caH(\frp\cap\frq)\otimes E\otimes {\caD'}\otimes S(\frk\cap \frh)\rightarrow S(\frg_+).$$
Let $\gamma'\colon Z(\frg_+)\rightarrow S(\frc_+)^{W(\frc_+)}$ be the
Harish-Chandra isomorphism.
Let ${\bf D'}=\bbC[D'_1,\ldots ,D'_r]$ be the polynomial subalgebra of $Z(\frg_+)$ generated by
$D'_i$ with $\gamma'(D'_i)=U'_i$   ($i=1,\cdots, r$).
Once again by comparing the grade with the filtration and arguing as in [W, 11.2.2], we
have the linear map
$$\Gamma_{\frg+}\colon \caH(\frp\cap\frq)\otimes E\otimes {\bf D'}\otimes U(\frk\cap\frh) \rightarrow U(\frg_+)$$
given by $ h^+\otimes e \otimes D' \otimes k\mapsto \symm(h^+)\symm(e)D'k$ is a linear isomorphism.

Let $(\tau,V_\tau)$ be an irreducible representation of $K$ and denote by $(\tau_+,V_\tau)$ its
restriction to $K_+=K\cap H$.
Let $\caE_\nu(G,\tau)$ (resp. $\caE_\nu(G_+,\tau_+)$) denote the space of
joint eigenfunctions of  $\bf D$ (resp. $\bf D'$) in the space of $\tau$-spherical functions $C^\infty(G,\tau)$
(resp. $C^\infty(G_+,\tau_+)$).
For each $f\in \caE_\nu(G,\tau)$, we define its Taylor series
$$T_f: U(\frg)\rightarrow V_\tau,  \text{ by }T_f(x)=xf(1), \text{ for }x\in U(\frg).$$
Let $\caH(\frp\cap \frq)$ be the harmonics for the symmetric algebra
$S(\frp\cap \frq)$ with the adjoint action of $K_+$ on $\frp\cap\frq$.
Hence, $T_f$ is completely determined by its evaluation on $\symm(\caH(\frp\cap\frq))\symm(E)$.
Set $$A_f(h^+\otimes e)=T_f(\symm(h^+)\symm(e)).$$
If $A_f=0$, then the function $f$ is constant zero.   Moreover,
$$A_f(ad(k)(h^+\otimes e))=
\tau_1(k)A_f(h^+\otimes e)\tau_2(k^{-1}), \text{ for }k\in K_+.$$  It follows that $A_f$ is contained in
$\Hom_{K_+}(\caH(\frp\cap\frq)\otimes E, V_\tau)$.
Similarly, for each $\phi\in \caE_\nu(G_+,\tau_+)$  the Taylor series $T_\phi$ is completely determined by its evaluation on $\symm(\caH(\frp\cap\frq))\symm(E)$.  In the proof of Theorem 4.3 in [HOW] it has been shown
that $\caE_\nu(G_+,\tau_+)$ is in linear bijection with $\Hom_{K_+}(\caH(\frp\cap\frq)\otimes E, V_\tau)$.
Define
$$B\colon \caE_\nu(G,\tau)\rightarrow \caE_\nu(G_+,\tau_+), \ f\mapsto \phi \ \text{iff}\ A_f=A_\phi.$$
Then $B$ is a linear injection.

\smallskip
Theorem 5.2. {\it $B$ is a linear bijection of  $\caE_\nu(G,\tau)$ with $\caE_\nu(G_+,\tau_+)$.}

\noindent {\it Proof.} It is enough to show $\dim \caE_\nu(G,\tau)= \dim \caE_\nu(G_+,\tau_+)$.
Note that
$$\dim \caE_\nu(G_+,\tau_+) =\dim \Hom_{K_+}(\caH(\frp\cap\frq)\otimes E, V_\tau)=
\dim E \cdot \dim \Hom_{K_+}(\caH(\frp\cap\frq), V_\tau),$$
which is equal to $\dim E \cdot \dim V_\tau^{M\cap K\cap H}$ and denoted by $d(\tau)$.
 We now show $\dim \caE_\nu(G,\tau)= d(\tau)$.  When $\nu$ is
generic,  the Eisenstein integrals corresponding to the matrix coefficients of irreducible principal series with $H$-fixed distribution vectors give $d(\tau)$ linear independent functions in $\caE_\nu(G,\tau)$.
Note that the $K$-types
of these principal series $\pi_{\xi,\lambda(\xi)_i}$ are copies of  $\Ind^K_M(\xi)$ (which does not depend on $\lambda(\xi)_i$,
$i=1,\dots,d=\dim E$) as following
$$\bigoplus_{i=1}^{d}\bigoplus_{w\in \caW}\bigoplus_{\xi}\Ind^K_M(\xi)=\bigoplus_{i=1}^d \bigoplus_{w\in \caW}\Ind^K_{w(K_M\cap H)w^{-1}}(1), $$
where  $\xi\in \widehat{M}_{fu}$ runs through discrete series of $M/wH_Mw^{-1}$.
Then the multiplicity of $V_\tau$ is
$$d\cdot \dim \Hom_K(\bigoplus_{w\in \caW}\Ind^K_{w(K_M\cap H)w^{-1}}(1),V_\tau)=d\cdot|\caW| \dim V_\tau^{M\cap K\cap H}=d(\tau).$$
Thus,  there exists
a basis $f_{1,\nu},\ldots,f_{d(\tau),\nu}$ for $\caE_\nu(G,\tau)$ such that each $f_{i,\nu}$ is holomorphic in generic $\nu$ and having
meromorphic extension to all $\nu$ ([B2]).  Let $D=\{z \in \bbC; |z|<1\}$ denote the unit disc
 and $D_0=\{z \in \bbC; 0<|z|<1\}$  the punctured disc. For $\nu_0$ non-generic, there exists a $\nu_1\in \fra^*_{\frq_\bbC}$ and such that $f_{i,z}=f_{i,\nu_0+z\nu_1}$ is holomorphic in $z\in D_0$ and  extended meromorphically to $D$. We apply Prop. 2.21 [OS1] to obtain linear independent $g_{i,z}=\sum_{j=1}^da_{ij}(z)f_{j,z}$ that are holomorphic in $z\in D$,
where $a_{ij}(z)$ are $d(\tau)^2$ meromorphic functions of $z\in D$. This extends the equality
$\dim \caE_\nu(G,\tau)= d(\tau)$ to $\nu=\nu_0$.  Thus, $B$ is a bijection.

This completes the proof of Theorems 5.2 as well as Theorem 5.4.

We remark that Theorems 7.5 and 8.4 can be proved with the above argument using a basis $f_{1,\nu},\ldots,f_{d(\tau),\nu}$ for  $\caE_\nu(G,\tau)$ with $\nu=\nu_0+z\nu_1$ and $f_{i,\nu_0+z\nu_1}$ being holomorphic in $z\in D$.  Still, we now prove the general statement in Theorem 5.5 on holomorphic dependence for the several-variables $\nu$.
\subsection*{4. A Proof of Theorem 5.5.}
Denote by $A$  the linear bijection defined above for the proof of Theorem 5.2, namely,
$$A\colon \caE_\nu(G,\tau)\rightarrow \Hom_{K_+}(\caH(\frp\cap\frq)\otimes E, V_\tau) \text{ by }  f\mapsto A_f. $$
Fix a basis $\eta_1\ldots,\eta_{d(\tau)}$ for $\Hom_{K_+}(\caH(\frp\cap\frq)\otimes E, V_\tau)$.
 Set $f_{i,\nu}=A^{-1}(\eta_i)$, namely $A_{f_{i,\nu}}=\eta_i$. We now prove that the basis $f_{1,\nu},\ldots,f_{d(\tau),\nu}$  for $\caE_\nu(G,\tau)$ satisfying the required condition, namely each $f_{i,\nu}$ is holomorphic in $\nu$.

First, we assume that $\nu_0$ is generic. There exists a basis $\{g_{i,\nu};i=1,\ldots ,d(\tau)\}$ of $\caE_\nu(G,\tau)$ such that each
$g_{i,\nu}$ is holomorphic in $\nu$ (in a domain $\Omega(\nu_0)$ of  $\nu_0$) as we argued above by using Eisenstein integrals.
Then $A({g_{i,\nu}})=\sum_jb_{ji}(\nu)\eta_j$ with $\nu\mapsto [b_{ij}(\nu)]$ a holomorphic map from $\Omega(\nu_0)$ to $GL(d(\tau),\bbC)$.
Let $[c_{ij}(\nu)]=[b_{ij}(\nu)]^{-1}$ be the inverse matrix.  Set $h_{i,\nu}=\sum_{j}c_{ji}(\nu)g_{j,\nu}$.
Then $h_{i,\nu}$ is holomorphic in $\nu$.  It follows that $h_{i,\nu}=f_{i,\nu}$,  since
 $A({h_{i,\nu}})=A({f_{i,\nu}})=\eta_i$.  Thus, $f_{i,\nu}$ is holomorphic at $\nu_0$.

If $\dim\fra^*_{\frq_\bbC}=1$, this has been done above by applying Prop. 2.21 [OS1].
Now we assume $\dim\fra^*_{\frq_\bbC}>1$.  We now show that
each $f_{i,\nu}$ is holomorphic in $\nu$ at any point $\nu_0\in \fra^*_{\frq_\bbC}$.
Note that the non-generic points $\caS$ is in the union of locally finite 
hypersurfaces $\cup \caH$ with each $\caH$ defined either by the zeros of the polynomial equations $J(\Lambda_\xi;x)=0$
or by  $\langle \nu,\alpha \rangle\in \Z$ for some restricted root $\alpha\in \Sigma(\frg,\fra_\frq)$.
Let $\caS_0\subseteq \caS$ be the points contained in precisely one of such hypersurfaces.
It follows from Corollary 7.3.2 \cite{K} that it suffices to show that
$f_{i,\nu}$ is holomorphic at any $\nu_0\in \caS_0$, since the singular locus of this extension is contained in $\caS\backslash\caS_0$ and the latter set has codimension
at least two in ${\fra_\frq}^*_\bbC$.

Suppose that $\nu_0\in \caS_0$ is in the hypersurface  $\caH=\{ \nu\in {\fra_\frq}^*_\bbC; g(\nu)=0\}$, where
the holomorphic function $g(\nu)$ is not constant zero. Recall that $D$  is the complex unit disc.
Then there is an open set $\Omega$ containing $\mu_0$ and a $\nu_1 \in \fra^*_{\frq_\bbC}$
so that the intersection of the disc $D(\mu)=\mu+\nu_1D$ with $\caS$ is contained in $\{\mu\}$ for any $\mu\in \Omega$.
Then there is an open set $\Omega(\nu_0)$ consisting of the union of discs
 $$\Omega(\nu_0)=\bigcup_{ \mu\in \caH\cap \Omega} D(\mu),$$
such that  $\Omega(\nu_0)\cap \caS\subset \caH$. Once again we  apply Prop. 2.21 [OS1] to obtain linear independent $g_{i,\nu}=\sum_{j=1}^{d(\tau)} a_{ij}(\nu)f_{j,\nu}$
that are holomorphic at any  $\nu\in \Omega(\nu_0)$,
where $a_{ij}(\nu)$ are $d(\tau)^2$ meromorphic functions of $\nu$ in $\Omega(\nu_0)$.  In particular, $g_{i,\nu}$ is
bounded at $\caH\cap \Omega(\nu_0)$.
By Riemmann's removable singularity theorem (Theorem 7.3.3 \cite{K}), $g_{i,\nu}$
is holomorphic at $\nu_0$.  Now we use  the same argument above for a generic point $\nu_0$ to get
$h_{i,\nu}$ holomorphic in $\nu\in \Omega(\nu_0)$ and $A({h_{i,\nu}})=\eta_i$.
Thus, $f_{i,\nu}=h_{i,\nu}$ is holomorphic at $\nu_0$.  This completes the proof of Theorem 5.5.

\subsection*{5. Correction to the list of exceptional simple pairs in Proposition 3.1.}
Recall that a simple symmetric pair $(\frg,\frh)$ is
called exceptional if the restriction $I(\frc)\rightarrow I(\fra_\frq)$ is not surjective.
There are totally 35 exceptional simple pairs which are listed in Proposition 3.1.
It follows from Helgason's restriction theorem [2] that a simple  pair  $(\frg,\frh)$ is exceptional if and only if
$$(\Sigma(\frg,\frc),\Sigma(\frg,\fra_\frq)) \in
\{(E_6,BC_2),(E_6,A_2),(E_7,C_3),(E_8,F_4)\}.$$
Denote by $(\frg^d,\frh^d)$ the dual symmetric pair of $(\frg,\frh)$.  Note that the two dual pairs have the
same pairs of the restriction root systems.
Thus, $(\frg,\frh)$ is exceptional if and only if $(\frg^d,\frh^d)$ is exceptional.   Since $G\times G/d(G)$ is dual to $G_\bbC/K_\bbC$, the following pairs
$$(\fre^\bbC_6,\fre_{6(-14)}), (\fre^\bbC_6,\fre_{6(-26)}), (\fre^\bbC_7,\fre_{7(-25)}), (\fre^\bbC_8,\fre_{8(-24)})$$
listed in Proposition 3.1 should be replaced by
$$(\fre^\bbC_6, \frs\fro_{10}(\bbC)+\bbC), (\fre^\bbC_6,\frf^\bbC_4), (\fre^\bbC_7,\fre^\bbC_6+\bbC), (\fre^\bbC_8,\fre^\bbC_7+\frs\frl_2(\bbC)).$$
The above four exceptional pairs are dual to $(\frg\times \frg,d(\frg))$
with $\frg=\fre_{6(-14)},\fre_{6(-26)},\fre_{7(-25)},\fre_{8(-24)}$.
Note that Riemannian symmetric pair $G/K$, $K_\epsilon$-symmetric space
$G/K_\epsilon$ and $G_\bbC/G_\bbR$ are self-dual [OS2].  The restriction for the 
exceptional pairs is calculated case by case in [HOW] Pages 640-642. 

\subsection*{6. A final remark.} Let ${\bbD}(G/H)$ be the algebra of $G$-invariant differential operators on an affine symmetric space $G/H$. The right action of $G$ on $C^\infty(G)$ induces a surjective
homomorphism
$$r\colon U(\frg)^H \rightarrow \bbD(G/H)$$
with kernel $U(\frg)^H\cap U(\frg)\frh$.  Note that $U(\frg)^H=(U(\frg)^H\cap U(\frg)\frh)\oplus \symm[S(\frq)^H]$.
 Then $r$ maps $\symm[S(\frq)^H]$ bijectively onto $\bbD(G/H)$.
Recall $\fra_\frq$ is a maximal abelian subspace in $\frp\cap \frq$
and  $\frb\supseteq \fra_\frq$ is a Cartan subspace of $G/H$.  Let
$$\gamma_\frb: \bbD(G/H)\rightarrow I(\frb)$$
be the Harish-Chandra isomorphism defined in [BS1].

A simple pair $(\frg,\frh)$ is called $\frb$-exceptional (in comparison with exceptional)
if the restriction $I(\frb)\rightarrow I(\fra_\frq)$ is not surjective.
It follows from Helgason's restriction theorem [2] that a simple  pair  $(\frg,\frh)$ is $\frb$-exceptional if and only if
$$(\Sigma(\frg,\frb),\Sigma(\frg,\fra_\frq)) \in
\{(E_6,BC_2),(E_6,A_2),(E_7,C_3),(E_8,F_4)\}.$$
There are totally 10 $\frb$-exceptional (in comparison with 35 exceptional) simple pairs:
$$(\fre_{6(-14)},\frs\frp_{2,2}), (\fre_{6(-26)},\frs\frp_{3,1}), (\fre_{7(-25)},\frs\fru_{6,2}),
(\fre_{7(-25)},\frs\fru^*_{8}), (\fre_{8(-24)},\frs\fro_{12,4}), (\fre_{8(-24)},\frs\fro^*_{16})$$
and four pairs of the form $(\frg\times \frg,d(\frg))$
with $\frg=\fre_{6(-14)},\fre_{6(-26)},\fre_{7(-25)},\fre_{8(-24)}$.

The main results of [1] depend on the choice of a polynomial algebra ${\bf D}\subseteq Z(\frg)$.
We can generalize the results to a polynomial algebra ${\bf D}\subset U(\frg)^H$ satisfying the following conditions:

(a) $\bf D\cong S(\fra_\frq)^{W(\fra_\frq)}$, if $(\frg,\frh)$ contains no simple
$\frb$-exceptional pair;

(b)  $\bf D$ is isomorphic to a subalgebra of $S(\fra_\frq)^{W(\fra_\frq)}$ such that $S(\fra_\frq)^{W(\fra_\frq)}$ is
a free $\bf D$-module  of finite rank, if $(\frg,\frh)$ contains some simple
$\frb$-exceptional pairs.

\noindent In particular, if $G/H$ is split, namely $\frb=\fra_\frq$, then ${\bf D}\cong {\bbD}(G/H)$.
The split symmetric spaces include (but not limited to) the following families:

(i) A Riemannian symmetric space $G/K$ and a $K_\epsilon$-symmetric space $G/K_\epsilon$;

(ii) A split groups $G$ regarded as a symmetric space $G\times G/d(G)$ and its dual space $G_\bbC/K_\bbC$;

(iii) The symmetric space $G_\bbC/ G_\bbR$ of a complex Lie group $G_\bbC$
over a real form  $G_\bbR$ that has a compact Cartan subgroup.

The details of the generalization will be published elsewhere.

\end{document}